\RequirePackage{ifpdf}
\ifpdf 
\documentclass[pdftex]{sigma}
\else
\documentclass{sigma}
\fi

\begin{document}

\allowdisplaybreaks
	
\renewcommand{\PaperNumber}{122}

\FirstPageHeading

\renewcommand{\thefootnote}{$\star$}

\ShortArticleName{Some Progress in Conformal Geometry}

\ArticleName{Some Progress in Conformal Geometry\footnote{This paper is a
contribution to the Proceedings of the 2007 Midwest
Geometry Conference in honor of Thomas~P.\ Branson. The full collection is available at
\href{http://www.emis.de/journals/SIGMA/MGC2007.html}{http://www.emis.de/journals/SIGMA/MGC2007.html}}}

\Author{Sun-Yung A. CHANG~$^\dag$, Jie QING~$^\ddag$ and Paul
YANG~$^\dag$}

\AuthorNameForHeading{S.-Y.A.~Chang, J.~Qing and P.~Yang}

\Address{$^\dag$~Department of Mathematics, Princeton University,
Princeton, NJ 08540, USA}
\EmailD{\href{mailto:chang@math.princeton.edu}{chang@math.princeton.edu}, \href{mailto:yang@math.princeton.edu}{yang@math.princeton.edu}}

\Address{$^\ddag$~Department of Mathematics, University of California, Santa Cruz,\\
$\phantom{^\ddag}$~Santa Cruz, CA
95064, USA}
\EmailD{\href{mailto:qing@ucsc.edu}{qing@ucsc.edu}}

\ArticleDates{Received August 30, 2007, in f\/inal form December
07, 2007; Published online December 17, 2007}

\Abstract{This is a survey paper of our current research on the
theory of partial dif\/feren\-tial equations in conformal geometry. Our
intention is to describe some of our current works in a rather brief
and expository fashion. We are not giving a comprehensive survey on
the subject and references cited here are not intended to be
complete. We introduce a~bubble tree structure to study the
degeneration of a class of Yamabe metrics on Bach f\/lat manifolds
satisfying some global conformal bounds on compact manifolds of
dimension 4. As applications, we establish a gap theorem, a
f\/initeness theorem for dif\/feomorphism type for this class, and
diameter bound of the $\sigma _2$-metrics in a class of conformal
4-manifolds. For conformally compact Einstein metrics we introduce
an eigenfunction compactif\/ication. As a consequence we obtain some
topological constraints in terms of renormalized volumes.}


\Classification{53A30; 53C20; 35J60}

\rightline{\it Dedicated to the memory of Thomas Branson}

\section[Conformal gap and finiteness theorem for a class of closed
4-manifolds]{Conformal gap and f\/initeness theorem\\ for a class of closed
4-manifolds}

\subsection{Introduction}

Suppose that $(M^4, g)$ is a closed 4-manifold. It follows from the
positive mass theorem that, for a 4-manifold with positive Yamabe
constant,
\[
\int_M \sigma_2 dv \leq 16\pi^2
\]
and equality holds if and only if $(M^4, g)$ is conformally
equivalent to the standard 4-sphere, where
\[
\sigma_2[g] = \frac 1{24}R^2 - \frac 12 |E|^2,
\]
$R$ is the scalar curvature of $g$ and $E$ is the traceless Ricci
curvature of $g$. This is an interesting fact in conformal geometry
because the above integral is a conformal invariant like the Yamabe
constant.

One may ask, whether there is a constant $\epsilon_0>0$ such that a
closed 4-manifold $M^4$ has to be dif\/feomorphic to $S^4$ if it
admits a metric $g$ with positive Yamabe constant and
\[
\int_M \sigma_2[g] dv_g \geq  (1-\epsilon)16\pi^2.
\]
for some $\epsilon < \epsilon_0$? Notice that here the Yamabe
invariant for such $[g]$ is automatically close to that for the
round 4-sphere. There is an analogous gap theorem of Bray and Neves
for Yamabe invariant in dimension 3 \cite{BN}. One cannot expect the
Yamabe invariant alone to isolate the sphere, and it is more
plausible to consider the integral of $\sigma_2$. We will answer the
question af\/f\/irmatively in the class of Bach f\/lat 4-manifolds.

Recall that Riemann curvature tensor decomposes into
\[
R_{ijkl} = W_{ijkl} + (A_{ik}g_{jl} - A_{il}g_{jk} - A_{jk}g_{il} +
A_{jl}g_{ik}),
\]
in dimension 4, where $\textsl{W}_{ijkl}$ is the Weyl curvature,
\[
A_{ij}= \frac 12 \left(R_{ij} - \frac 16 Rg_{ij}\right)
\]
is Weyl--Schouten curvature tensor and $R_{ij}$ is the Ricci
curvature tensor. Also recall that the Bach tensor is
\[
B_{ij}= W_{kijl,lk} + \frac 12 R_{kl} W_{kijl}.
\]
We say that a metric $g$ is Bach f\/lat if $B_{ij}=0$. Bach f\/lat
metrics are critical metrics for the functional $\int_M |W|^2 dv$.
Bach f\/latness is conformally invariant in dimension 4. It follows
from Chern--Gauss--Bonnet,
\[
8\pi^2 \chi(M^4)= \int_M (\sigma_2 + |W|^2)dv,
\]
that $\int_M \sigma_2 dv$ is conformally invariant.

The gap theorem is as follows:

\begin{theorem}\label{theorem 1.1.1} Suppose that $(M^4, [g])$ is a Bach flat
closed $4$-manifold with positive Yamabe constant and that
\[
\int_M (|W|^2dv)[g] \leq \Lambda_0
\]
for some fixed positive number $\Lambda_0$. Then there is a positive
number $\epsilon_0 > 0$ such that, if
\[
\int_M \sigma_2[g] dv_g \geq (1-\epsilon)16\pi^2
\]
holds for some constant $\epsilon < \epsilon_0$, then $(M^4, [g])$
is conformally equivalent to the standard $4$-sphere.
\end{theorem}

Our approach is based on the recent work on the compactness of Bach
f\/lat metrics on 4-manifolds of Tian and Viaclovsky~\cite{TV-1,
TV-2}, and of Anderson~\cite{An}. Indeed our work relies on a more
precise understanding of the bubbling process near points of
curvature concentration. For that purpose we develop the bubble tree
structure in a sequence of metrics that describes precisely the
concentration of curvature. Our method to develop bubble tree
structure is inspired by the work of Anderson and Cheeger \cite{AC}
on the bubble tree conf\/igurations of the degenerations of metrics of
bounded Ricci curvature. Our construction is modeled after this work
but dif\/fers in the way that our bubble tree is built from the
bubbles at points with the smallest scale of concentration to
bubbles with larger scale; while the bubble tree in \cite{AC} is
constructed from bubbles of large scale to bubbles with smaller
scales. The inductive method of construction of our bubble tree is
modeled on earlier work of \cite{BC, Q, St} on the study of
concentrations of energies in harmonic maps and the scalar curvature
equations.

As a consequence of the bubble tree construction we are able to
obtain the following f\/inite dif\/feomorphism theorem:

\begin{theorem}\label{theorem 1.1.2} Suppose that $\textbf{A}$ is a collection of Bach
flat Riemannian manifolds $(M^4, g)$ with positive Yamabe constant,
satisfying
\[
\int_M (|W|^2 dv)[g] \leq \Lambda_0,
\]
for some fixed positive number $\Lambda_0$, and
\[
\int_M (\sigma_2 dv)[g] \geq \sigma_0,
\]
for some fixed positive number $\sigma_0$. Then there are only
finite many diffeomorphism types in $\textbf{A}$.
\end{theorem}

It is known that in each conformal class of metrics belonging to the
family $\textbf{A}$, there is a metric $\bar g=e^{2w}g$ such that
$\sigma _2(A_{\bar g}) =1$, which we shall call the $\sigma _2$
metric. The bubble tree structure in the degeneration of Yamabe
metrics in $\textbf{A}$ is also helpful to understand the behavior
of the $\sigma_2$-metrics in $\textbf{A}$. For example:

\begin{theorem}\label{theorem 1.1.3} For the conformal classes $\left[g_0\right]
\in {\textbf{A}}$ the conformal metrics $g=e^{2w}g_0$ satisfying the
equation $\sigma _2(g)=1$ has a uniform bound for the diameter.
\end{theorem}

The detailed version of this work has appeared in our paper
\cite{CQY1}.

\subsection{The neck theorem}

The main tool we need to develop the bubble tree picture is the
 neck theorem  which should be compared with the neck
theorem in the work of Anderson and Cheeger \cite{AC}. Due to the
lack of point-wise bounds on Ricci curvature, our version of the
neck theorem will have weaker conclusion. But it is suf\/f\/icient to
allow us to construct the bubble tree at each point of curvature
concentration.

Let $(M^4, g)$ be a Riemannian manifold. For a point $p\in M$,
denote by $B_r(p)$ the geodesic ball with radius $r$ centered at
$p$, $S_r(p)$ the geodesic sphere of radius $r$ centered at $p$.
Consider the geodesic annulus centered at $p$:
\[
\bar A_{r_1, r_2} (p) = \{ q\in M: r_1 \leq \text{dist}(q, p) \leq
r_2\}.
\]
In general, $\bar A_{r_1, r_2}(p)$ may have more than one connected
components. We will consider any one component
\[
A_{r_1,r_2}(p) \subset \bar A{r_1, r_2}(p)
\]
that meets the geodesic sphere of radius $r_2$:
\[
A_{r_1,r_2}(p)\bigcap S_{r_2}(p)\neq \varnothing.
\]
Let $H^3(S_r(p))$ be the 3D-Hausdorf\/f measure of the geodesic sphere
$S_r(p)$.

\begin{theorem}\label{theorem 1.2.1} Suppose $(M^4, g)$ is a Bach flat and
simply connected $4$-manifold with a Yamabe metric of positive Yamabe
constant. Let $p \in M$, $\alpha \in (0,1)$, $\epsilon
>0$, $v_1 >0$, and $a < \text{\rm dist}(p, \partial M)$.
Then there exist positive numbers $\delta_0$, $c_2$, $n$ depending on
$\epsilon$, $\alpha$, $C_s$, $v_1$, $a$ such that the following holds. Let
$A_{r_1, r_2}(p)$ be a connected component of the geodesic annulus
in $M$ such that
\begin{gather*}
r_2 \leq c_2 a, \qquad r_1 \leq \delta_0 r_2,
\\
H^3 (S_r(p)) \leq v_1 r^3, \quad\forall \; r\in [r_1, 100r_1],
\end{gather*}
and
\[
\int_{A_{r_1, r_2}(p)} |Rm|^2 dv \leq \delta_0.
\]
Then $A_{r_1, r_2}(p)$ is the only such component. In addition,
 for the only component
\[
A_{(\delta_0^{-\frac 14}-\epsilon)r_1, (\delta_0^\frac 14 +\epsilon)
r_2}(p),
\]
which intersects with $S_{(\delta_0^\frac 14 +\epsilon) r_2}(p)$,
there exist some $\Gamma \subset O(4)$, acting freely on $S^3$, with
$|\Gamma| \leq n$, and an quasi isometry $\psi$, with
\[
A_{(\delta_0^{-\frac 14}+\epsilon)r_1, (\delta_0^\frac 14 -\epsilon)
r_2}(p)
\subset \Psi(C_{\delta_0^{-\frac 14}r_1,
\delta_0^\frac 14 r_2} (S^3/\Gamma))
\subset
A_{(\delta_0^{-\frac 14}-\epsilon)r_1, (\delta_0^\frac 14 +\epsilon)
r_2}(p)
\]
such that for all $C_{\frac 12 r, r}(S^3/\Gamma) \subset
C_{\delta_0^{-\frac 14}r_1, \delta_0^\frac 14 r_2} (S^3/\Gamma)$, in
the cone  coordinates, one has
\[
|(\Psi^*(r^{-2}g))_{ij} - \delta_{ij}|_{C^{1, \alpha}} \leq
\epsilon.
\]
\end{theorem}

The f\/irst step in the proof is to use the Sobolev inequality to show
the uniqueness of the connected annulus $A_{r_1, r_2}(p)$. The
second step is to establish the growth of volume of geodesic spheres
\[
H^3(S_r(p) \leq C r^3
\]
for all $r\in [r_1, \frac 12 r_2]$. Here we rely on the work of Tian
and Viaclovsky \cite{TV-1, TV-2} where they analyzed the end
structure of a Bach-f\/lat, scalar f\/lat manifolds with f\/inite $L^2$
total curvature. The last step is to use the Gromov and Cheeger
compactness argument as in the work of Anderson and Cheeger~\cite{AC} to get the cone structure of the neck.

\subsection{Bubble tree construction}

In this section we attempt to give a clear picture about what
happen at curvature concentration points. We will detect and extract
bubbles by locating the centers and scales of curvature concentration.

We will assume here that $(M_i, g_i)$ are Bach f\/lat 4-manifolds with
positive scalar curvature Yamabe metrics, vanishing f\/irst homology,
and f\/inite $L^2$ total curvature.  Choose $\delta$ small enough
according to the $\epsilon$-estimates and the neck theorem in the
previous section. Suppose that $X_i \subset M_i$ contains a geodesic
ball of a f\/ixed radius $r_0$ and
\[
\int_{T_{\eta_0}(\partial X_i)} |Rm^i|^2 dv^i \leq \frac \delta 2,
\]
where
\[
T_{\eta_0}(\partial X_i) = \{p\in M_i: \text{dist}(p,
\partial X_i) < \eta_0 \},
\]
for some f\/ixed positive number $4\eta_0 < r_0$. Def\/ine, for $p\in
X_i$,
\[
s_i^1(p) = r \qquad \text{such that} \quad \int_{B^i_r(p)} |Rm^i|^2dv^i =
\frac \delta 2.
\]
Let
\[
p_i^1 = p  \qquad \text{such that} \quad s^1_i(p) = \inf_{B^i_{t_0}(p_i)}
s_i^1(p).
\]
We may assume $\lambda_i^1 = s_i^1(p_i^1)\to 0$, for otherwise there would be no
curvature concentration in $X_i$. We then conclude that $(M_i,
(\lambda_i^1)^{-2}g_i, p_i^1)$ converges to $(M_\infty^1,
g_\infty^1, p_\infty^1)$, which is a Bach f\/lat, scalar f\/lat,
complete 4-manifold satisfying the Sobolev inequality, having f\/inite $L^2$ total
curvature, and one single end.

\begin{definition}\label{definition 1.3.1} We call a Bach f\/lat,
scalar f\/lat, complete 4-manifold with the
Sobolev inequality, f\/inite $L^2$ total curvature, and a single ALE
end a {\it leaf bubble}, while we will call such space with f\/initely
many isolated irreducible orbifold points an {\it intermediate
bubble}.
\end{definition}

Now, we def\/ine, for $p\in X_i \setminus B^i_{K^1\lambda_i^1}
(p_i^1)$,
\[
s_i^2(p) = r
\]
such that
\[
\int_{B^i_r (p) \setminus B^i_{K^1\lambda_i^1} (p_i^1)} |Rm^i|^2dv^i
= \frac \delta 2.
\]
Let
\[
p_i^2 = p
\]
such that
\[
s^2_i(p) = \inf_{B^i_r (p)\setminus B^i_{K^1\lambda_i^1} (p_i^1)}
s_i^2(p).
\]
Again let $\lambda_i^2 = s_i^2(p_i^2)\to 0$. Otherwise there would
be no more curvature concentration. Then

\begin{lemma}\label{lemma 1.3.1}
\[
\frac {\lambda_i^2}{\lambda_i^1} + \frac {\text{\rm dist}(p_i^1,
p_i^2)}{\lambda_i^1} \to \infty.
\]
\end{lemma}

There are two possibilities:
\begin{gather*}
 \text{Case 1.} \qquad  \frac {\text{dist}(p_i^1,
p_i^2)}{\lambda_i^2} \to \infty; \\
 \text{Case 2.} \qquad \frac {\text{dist}(p_i^1,
p_i^2)}{\lambda_i^2} \leq M^1.
\end{gather*}
In Case 1, we certainly also have
\[
\frac {\text{dist}(p_i^1, p_i^2)}{\lambda_i^1} \to \infty.
\]
Therefore, in the convergence of the sequence $(M_i, (\lambda_i^2)^{-2}g_i, p_i^2)$
the concentration which produces the bubble $(M_\infty^1,
g_\infty^1)$ eventually escapes to inf\/inity of $M_\infty^2$
and hence is not visible to the
bubble $(M_\infty^2, g_\infty^2)$, likewise, in the converging sequence
$(M_i, (\lambda_i^1)^{-2}g_i, p_i^1)$ one does
not see the concentration which produces $(M_\infty^2, g_\infty^2)$.
There are at most f\/inite number of such leaf bubbles.

\begin{definition}\label{definition 1.3.2} We say two bubbles $(M_\infty^{j_1},
g_\infty^{j_1})$ and $(M_\infty^{j_2}, g_\infty^{j_2})$ associated
with $(p_i^{j_1}, \lambda_i^{j_1})$ and $(p_i^{j_2},
\lambda_i^{j_2})$ are {\it separable} if
\[
\frac{\text{dist}(p_i^{j_1}, p_i^{j_2})}{\lambda_i^{j_1}} \to \infty
\qquad \text{and} \qquad \frac{\text{dist}(p_i^{j_1},
p_i^{j_2})}{\lambda_i^{j_2}} \to \infty.
\]
\end{definition}

In Case 2, one starts to trace intermediate bubbles which will be
called parents of some bubbles. We would like to emphasize a very
important point here. One needs the neck theorem to take limit in
Goromov--Hausdorf\/f topology to produce the intermediate bubbles. The
neck Theorem is used to prove the limit space has only isolated
point singularities, which are then proven to be orbifold points.

\begin{lemma}\label{lemma 1.3.2} Suppose that there are several separable bubbles
$\{(M_\infty^j, g_\infty^j)\}_{j\in J}$ associated with $\{(p_i^j,
\lambda_i^j)\}_{j\in J}$. Suppose that there is a concentration
detected as $(p_i^k, \lambda_i^k)$ after $\{(p_i^j,
\lambda_i^j)\}_{j\in J}$  such that
\[
\frac {\text{\rm dist}(p_i^k, p_i^j)}{\lambda_i^k} \leq M^j,
\]
therefore
\[
\frac {\lambda_i^k}{\lambda_i^j} \to \infty
\]
for each $j\in J$. In addition, suppose that $\{(p_i^j, \lambda_i^j)\}_{j\in
J}$ is the maximal collection of such. Then
\newline $(M_i, (\lambda_i^k)^{-2}g_i, p_i^k)$ converges in
Gromov--Hausdorff topology to an intermediate bubble $(M_\infty^k,
g_\infty^k)$. $(M_\infty^k, g_\infty^k)$ is either a parent or a
grandparent of all the given bubbles $\{(M_\infty^j,
g_\infty^j)\}_{j\in J}$.
\end{lemma}

We remark that it is necessary to create some strange intermediate
bubbles to handle the inseparable bubbles. This situation does not
arise in the degeneration of Einstein metrics. In that case there is
a gap theorem for Ricci f\/lat complete orbifolds and there is no
curvature concentration at the smooth points due to a simple volume
comparison argument, both of which are not yet available in our
current situation. We will call those intermediate bubbles {\it
exotic bubbles}.

\begin{definition}\label{definition 1.3.3} A {\it bubble tree} $T$ is def\/ined to be a
tree whose vertices are bubbles and whose edges are necks from neck
Theorem. At each vertex $(M_\infty^j, g_\infty^j)$, its ALE end is
connected, via a~neck, to its parent towards the {\it root bubble}
of $T$, while at f\/initely many isolated possible orbifold points of
$(M_\infty^j, g_\infty^j)$, it is connected, via necks, to its
children towards leaf bubbles of $T$. We say two bubble trees $T_1$
and $T_2$ are {\it separable} if their root bubbles are separable.
\end{definition}

To f\/inish this process we just iterate the process of extracting
bubbles the construction has to end at some f\/inite steps. In summary
we have

\begin{theorem}\label{theorem 1.3.2} Suppose that $(M_i, g_i)$ are Bach flat
$4$-manifolds with positive scalar curvature Yamabe metrics, vanishing
first homology, and finite $L^2$ total curvature. Then $(M_i, g_i)$
converges to Bach-flat $4$-manifold $(M_\infty, g_\infty)$ with
finitely orbifold singularities $S$. The convergence is strong in
$C^\infty$ away from a finite number of points $B\supset S$. At each
point $b$ in $B$ there is a bubble tree attached to $b$.
\end{theorem}

\section{Conformally compact Einstein manifolds}

\subsection{Conformally compact Einstein manifolds}

Suppose that $X^{n+1}$ is a smooth manifold of dimension $n+1$ with
smooth boundary $\partial X = M^n$. A def\/ining function for the
boundary $M^n$ in $X^{n+1}$ is a smooth function $x$ on $\bar
X^{n+1}$ such that
\[
\left\{\begin{array}{ll} x>0 \ \ & \text{in $X$};  \\
x=0 \ \ & \text{on $M$};  \\
dx \neq 0 \ \ & \text{on $M$}.
\end{array}\right.
\]
A Riemannian metric $g$ on $X^{n+1}$ is conformally compact if
$(\bar X^{n+1},  x^2 g)$ is a compact Riemannian manifold with
boundary $M^n$ for a def\/ining function $x$. Conformally compact
manifold $(X^{n+1},  g)$ carries a well-def\/ined conformal structure
on the boundary $M^n$, where each metric $\hat g$ in the class is
the restriction of $\bar g = x^2g$ to the boundary $M^n$ for a
def\/ining function $x$. We call $(M^n, [\hat g])$ the conformal
inf\/inity of the conformally compact manifold $(X^{n+1}, g)$. A short
computation yields that, given a def\/ining function $x$,
\[
R_{ijkl}[g] = |dx|^2_{\bar g} (g_{ik}g_{jl}-g_{il}g_{jk}) + O(x^3)
\]
in a coordinate $(0, \epsilon)\times M^n \subset X^{n+1}$.
Therefore, if we assume that $g$ is also asymptotically hyperbolic,
then
\[
|dx|^2_{\bar g}|_M =1
\]
for any def\/ining function $x$. If $(X^{n+1},  g)$ is a conformally
compact manifold and $\text{Ric}[g] = -ng$, then we call $(X^{n+1}, g)$ a conformally compact Einstein manifold.

Given a conformally compact, asymptotically hyperbolic manifold
$(X^{n+1}, g)$ and a representative $\hat g$ in $[\hat g]$ on the
conformal inf\/inity $M^n$, there is a uniquely determined def\/ining function $x$
such that, on $M\times (0, \epsilon)$ in $X$, $g$ has the normal
form
\begin{equation}
g = x^{-2} (dx^2 + g_x), \label{equation1}
\end{equation}
where $g_x$ is a 1-parameter family of metrics on $M$. This is
because

\begin{lemma}\label{lemma 2.1.1} Suppose that $(X^{n+1}, g)$ is a conformally
compact, asymptotically hyperbolic manifold with the conformal
infinity $(M, [\hat g])$. Then, for any $\hat g \in [\hat g]$, there
exists a unique defining function~$x$ such that
\[
|dx|^2_{r^2g} = 1
\]
in a neighborhood of the boundary $[0, \epsilon)\times M$ and
\[
r^2 g |_M = \hat g.
\]
\end{lemma}

Given a conformally compact Einstein manifold $(X^{n+1}, g)$, in the
local product coordinates $(0, \epsilon)\times M^n$ near the
boundary where the metric takes the normal form \eqref{equation1}, the Einstein
equations split and display some similarity to a second order
ordinary dif\/ferential equations with a regular singular point.

\begin{lemma}\label{lemma 2.1.2} Suppose that $(X^{n+1}, g)$ is a conformally
compact Einstein manifold with the conformal infinity $(M^n, [\hat
g])$ and that $x$ is the defining function associated with a metric
$\hat g\in [\hat g]$. Then
\begin{gather*}
 g_x  = \hat g + g^{(2)}x^2 + (\text{even powers of $x$})
 + g^{(n-1)}x^{n-1} + g^{(n)}x^n + \cdots,
\end{gather*}
when $n$ is odd, and
\begin{gather*} g_x  = \hat g + g^{(2)}x^2 + (\text{even powers of $x$})
 + g^{(n)}x^n + h x^n\log x + \cdots,
\end{gather*}
when $n$ is even, where:

a) $g^{(2i)}$ are determined by $\hat g$ for $2i < n$;

b) $g^{(n)}$ is traceless when $n$ is odd;

c) the trace part of $g^{(n)}$ is determined by $\hat g$ and $h$ is
traceless and determined by $\hat g$;

d) the traceless part of $g^{(n)}$ is divergence free.
\end{lemma}

Readers are referred to \cite{G} for more details about the above
two lemmas.

\subsection{Examples of conformally compact Einstein
manifolds}\label{sec2.2}

Let us look at some examples.

\medskip

\noindent{\it a) The hyperbolic spaces}
\[
\left(R^{n+1},  \frac {(d|x|)^2}{1+|x|^2} + |x|^2 d\sigma\right),
\]
where $d\sigma$ is the standard metric on the $n$-sphere. We may
write
\[
g_H = s^{-2}\left(ds^2 + \left(1-\frac {s^2}4\right)^2 d\sigma\right),
\]
where
\[
s = \frac 2{\sqrt{1+|x|^2} + |x|}
\]
is a def\/ining function. Hence the conformal inf\/inity is the standard
round sphere $(S^n, d\sigma)$.

\medskip

\noindent{\it b) The hyperbolic manifolds}
\[
\left(S^1 (\lambda)\times R^n,  (1+r^2)dt^2 + \frac {dr^2}{1+r^2} +
r^2d\sigma\right).
\]
Let
\[
r = \frac {1-\frac {s^2}4}s = \sinh \log \frac 2s
\]
for a def\/ining function $s$. Then
\[
g_H^0 = s^{-2}\left(ds^2 + \left(1+\frac {s^2}4\right)^2dt^2+ \left(1-\frac
{s^2}4\right)^2d\sigma \right).
\]
Thus the conformal inf\/inity is standard $(S^1(\lambda)\times
S^{n-1}, dt^2 + d\sigma)$.

\medskip
 \noindent{\it c) AdS-Schwarzchild}
\[
\big(R^2\times S^2,  g_{+1}^m\big),
\]
where
\[
g_{+1}^m = V dt^2 + V^{-1}dr^2 + r^2 g_{S^2},
\qquad
V = 1+r^2 - \frac {2m}r,
\]
$m$ is any positive number, $ r\in [r_h, +\infty)$, $t\in
S^1(\lambda)$ and $(\theta, \phi)\in S^2$, and $r_h$ is the positive
root for $1 + r^2 - \frac {2m}r =0$. In order for the metric to be
smooth at each point where $S^1$ collapses we need  $Vdt^2 +
V^{-1}dr^2$ to be smooth at $r= r_h$, i.e.
\[
V^\frac 12 \frac {d(V^\frac 12 2\pi\lambda)}{dr} \Big|_{r=r_h} = 2\pi.
\]
Note that its conformal inf\/inity is $(S^1(\lambda)\times S^2, [dt^2
+ d\theta^2 + \sin^\theta d\phi^2])$ and $S^1$ collapses at the
totally geodesic $S^2$, which is the so-called horizon.
Interestingly, $\lambda$ is does not vary monotonically  in $r_h$,
while $r_h$ monotonically depends on $m$. In fact, for each $0
<\lambda < 1/\sqrt 3$, there are two dif\/ferent~$m_1$ and~$m_2$ which
share the same $\lambda$. Thus, for the same conformal inf\/inity
$S^1(\lambda)\times S^2$ when $0 < \lambda < 1/\sqrt 3$, there are
two non-isometric AdS-Schwarzschild space with metric $g^+_{m_1}$
and~$g^+_{m_2}$ on $R^2\times S^2$. These are the interesting simple
examples of non-uniqueness for conformally compact Einstein metrics.

\medskip
\noindent{\it d) AdS-Kerr spaces}
\[
\big(CP^2\setminus \{p\}, g_\alpha\big),
\]
where $p$ is a point on $CP^2$,
\begin{gather*}
 g_\alpha  = E_\alpha((r^2-1)F_\alpha^{-1}dr^2 +
(r^2-1)^{-1}F_\alpha (dt + \cos \theta d\phi)^2
 + (r^2-1)(d\theta^2 + \sin^2\theta d\phi^2)),
\\
E_\alpha = \frac 23 \frac {\alpha-2}{\alpha^2 -1},
\qquad
F_\alpha = (r-\alpha)((r^3 - 6r + 3\alpha^{-1})E_\alpha +
4(r-\alpha^{-1})),
\end{gather*}
$r\geq \alpha$, $t\in S^1(\lambda)$, and $(\theta, \phi)\in S^2$.
For the metric to be smooth at the horizon, the totally geodesic
$S^2$, we need to require
\[
\sqrt{\frac{F}{E(r^2-1)}}\frac{d}{dr}\left(2\pi \lambda\sqrt{\frac
{EF}{r^2-1}}\right) = 2\pi.
\]
Here $(t,\theta, \phi)$ is the coordinates for $S^3$ through the
Hopf f\/iberation. The conformal inf\/inity is the Berger sphere with
the Hopf f\/ibre of length $\pi E_\alpha$ and the $S^2$ of area $4\pi
E_\alpha$. For every $0 < \lambda < (2-\sqrt 3)/3$ there are exactly
two $\alpha$, hence two AdS-Kerr metrics $g_\alpha$. It is interesting
to note that
$(2-\sqrt 3)/3 < 1$, so the standard $S^3(1)$ is not included in
this family.

One may ask, given a conformal manifold $(M^n,  [\hat g])$, is
there a conformally compact Einstein manifold $(X^{n+1},  g)$ such
that $(M^n,  [\hat g])$ is the conformal inf\/inity? This in general
is a dif\/f\/icult open problem. Graham and Lee in \cite{GL} showed that
for any conformal structure that is a~perturbation of the round one
on the sphere $S^n$ there exists a conformally compact Einstein
metric on the ball $B^{n+1}$.

\subsection[Conformal compactifications]{Conformal compactif\/ications}

Given a conformally compact Einstein manifold $(X^{n+1}, g)$, what
is a good conformal compactif\/ication? Let us consider the hyperbolic
space. The hyperbolic space $(H^{n+1}, g_H)$ is the hyperboloid
\[
\big\{(t, x)\in R\times R^{n+1}: -t^2 + |x|^2 = -1, t> 0\big\}
\]
in the Minkowski space-time $R^{1, n+1}$. The stereographic
projection via the imaginary south pole gives the Poincar\'{e} ball
model
\[
\left(B^{n+1}, \left(\frac 2{1-|y|^2}\right)^2|dy|^2\right)
\]
and replacing the $x$-hyperplane by $z$-hyperplane tangent to the
light cone gives the half-space model
\[
\left(R^{n+1}_+, \frac {|dz|^2}{z_{n+1}^2}\right),
\]
where
\[
\frac {1+|y|^2}{1-|y|^2} = t, \qquad \frac 1{z_{n+1}} = t -
x_{n+1}.
\]
Therefore
\begin{gather*}
(H^{n+1}, t^{-2}g_H)=(S^{n+1}_+, g_{S^{n+1}}),
\\
(H^{n+1}, (t+1)^{-2}g_H)=(B^{n+1}, |dy|^2),
\\
(H^{n+1}, (t-x_{n+1})^{-2}g_H) = (R^{n+1}_+, |dz|^2).
\end{gather*}
The interesting fact here is that all coordinate functions $\{t,
x_1, x_2, \dots, x_{n+1}\}$ of the Minkowski space-time are
eigenfunctions on the hyperboloid. Thus positive eigenfunctions on a
conformally compact Einstein manifold are expected to be candidates
for good conformal compactif\/ications. This is f\/irst observed in
\cite{Q1}.

\begin{lemma}\label{lemma 2.3.1} Suppose that $(X^{n+1}, g)$ is a conformally
compact Einstein manifold and that $x$ is a special defining
function associated with a representative $\hat g\in [\hat g]$. Then
there always exists a~unique positive eigenfunction $u$ such that
\[
\Delta u = (n+1) u \qquad \text{in} \ X
\]
and
\[
u = \frac 1x + \frac {R[\hat g]}{4n(n-1)}x + O (x^2)
\]
near the infinity.
\end{lemma}

We remark here that, for the hyperbolic space $H^{n+1}$ and the
standard round metric in the inf\/inity, we have
\[
t = \frac 1x + \frac 14 x.
\]
As we expect, positive eigenfunctions indeed give a preferable
conformal compactif\/ication.

\begin{theorem}\label{theorem 2.3.2} Suppose that $(X^{n+1},g)$ is a conformally
compact Einstein manifold, and that $u$ is the eigenfunction
obtained for a Yamabe metric $\hat g$ of the conformal infinity $(M,
[\hat g])$ in the previous lemma. Then $(X^{n+1}, u^{-2}g)$ is a
compact manifold with totally geodesic boundary $M$ and
\[
R[u^{-2} g] \geq  \frac {n+1}{n-1} R[\hat g].
\]
\end{theorem}

As a consequence

\begin{corollary}\label{corollary 2.3.3} Suppose that $(X^{n+1}, g)$ is a conformally
compact Einstein manifold and its conformal infinity is of positive
Yamabe constant. Suppose that $u$ is the positive eigenfunction
associa\-ted with the Yamabe metric on the conformal infinity obtained
in Lemma~{\rm \ref{lemma 1.3.1}}. Then $(X^{n+1}, u^{-2}g)$ is a compact manifold with
positive scalar curvature and totally geodesic boundary.
\end{corollary}

The work of Schoen--Yau and Gromov--Lawson then give some topological obstruction for a
conformally compact Einstein manifold to have its conformal inf\/inity
of positive Yamabe constant.
A surprising consequence of the eigenfunction compactif\/ications is the
rigidity of the hyperbolic space without assuming the spin
structure.

\begin{theorem}\label{theorem 2.3.4} Suppose that $(X^{n+1}, g)$ is a conformally
compact Einstein manifold with the round sphere as its conformal
infinity. Then $(X^{n+1}, g)$ is isometric to the hyperbolic space.
\end{theorem}

\subsection{Renormalized volume}

We will introduce the renormalized volume, which was f\/irst
noticed by physicists in their investigations of the holography
principles in AdS/CFT. Take a def\/ining function $x$ associated with
a choice of the metric $\hat g\in [\hat g]$ on the conformal
inf\/inity, then compute, when $n$ is odd,
\begin{gather}\label{4.1}
\text{Vol}(\{ x > \epsilon\} = c_0 \epsilon^{-n} +  \text{odd
powers of $\epsilon$}    + V + o(1),
\end{gather}
when $n$ is even,
\begin{gather}\label{4.2}
\text{Vol}(\{ s > \epsilon\}) = c_0 \epsilon^{-n} +  \text{even
powers of $\epsilon$} + L\log \frac 1\epsilon + V + o(1).
\end{gather}
It turns out the numbers $V$ in odd dimension and $L$ in even dimension are
independent of the choice of the metrics in the class. We will see
that $V$ in even dimension is in fact a conformal anomaly.

\begin{lemma}\label{lemma 2.4.1} Suppose that $(X^{n+1}, g)$ is a conformally
compact Einstein manifold and that $\bar x$ and~$x$ are two defining
functions associated with two representatives in~$[\hat g]$ on the
conformal infinity $(M^n, [\hat g])$. Then
\[
\bar x = x e^w
\]
for a function $w$ on a neighborhood of the boundary $[0,
\epsilon)\times M$ whose expansion at $x=0$ consists of only even
powers of $x$ up through and including $x^{n+1}$ term.
\end{lemma}

\begin{theorem}\label{theorem 2.4.2} Suppose that $(X^{n+1}, g)$ is a conformally
compact Einstein manifold. The $V$ in~\eqref{4.1} when $n$ is odd and $L$
in~\eqref{4.2} when n$n$ is even are independent of the choice of
representative $\hat g\in [\hat g]$ on the conformal infinity $(M^n,
[\hat g])$.
\end{theorem}

Let us calculate the renormalized volume for the examples in Section~\ref{sec2.2}.

\medskip

\noindent{\it a) The hyperbolic space}: \quad We recall
\[
(H^4, g_H) = \left(B^4, \left(\frac 2{1-|y|^2}\right)^2 |dy|^2\right),
\]
where
\[
g_H = s^{-2}\left(ds^2 + \left(1-\frac {s^2}4\right)^2 h_0\right)
\]
and $h_0$ is the round metric on $S^3$. Then
\[
\text{vol}(\{ s >\epsilon\}) = \int_\epsilon^2 \int_{S^3}s^{-4}
\left(1-\frac {s^2}4\right)^3d\sigma_0ds
\]
where $d\sigma_0$ is the volume element for the round unit sphere
\begin{gather*}
 \text{vol}(\{ s >\epsilon\})  = 2\pi^2 \int_\epsilon^2
s^{-4}\left( 1 - \frac {3s^2}4 + \frac {3s^4}{16} - \frac {s^6}{64}\right)ds
\\
\phantom{\text{vol}(\{ s >\epsilon\})}{}  = 2\pi^2\left( -\frac 13 s^{-3}\Big|_\epsilon^2 + \frac 34
s^{-1}\Big|_\epsilon^2 + \frac 3{16}(2-\epsilon) - \frac 1{3\times
64}s^3\Big|_\epsilon^2\right) \\
\phantom{\text{vol}(\{ s >\epsilon\})}{} = \frac {2\pi^2}3 \epsilon^{-3} - \frac {3\pi^2}2 \epsilon^{-1} +
2\pi^2\left(- \frac 1{3\times 8} + \frac 38 + \frac 38 - \frac
1{3\times 8}\right) + O(\epsilon) \\
\phantom{\text{vol}(\{ s >\epsilon\})}{} = \frac {2\pi^2}3 \epsilon^{-3} - \frac {3\pi^2}2 \epsilon^{-1} +
\frac {4\pi^2}3 + O(\epsilon).
\end{gather*}
Thus
\[
V(H^4, g_H) = \frac {4\pi^2}3.
\]

\medskip
\noindent{\it b) The hyperbolic manifold}: \quad We
recall
\[
\left(S^1 (\lambda)\times R^3, (1+r^2)dt^2 + \frac {dr^2}{1+r^2} +
r^2g_{S^2}\right)
\]
and
\[
g_H^0 = s^{-2}\left(ds^2 + \left(1-\frac {s^2}4\right)^2(d\theta^2 + \sin^2\theta
d\phi^2) + \left(1+\frac {s^2}4\right)^2dt^2\right).
\]
Then
\[
\text{vol}(\{ s > \epsilon\}) = \int_\epsilon^2 \int_{S^2}\int_{S^1}
s^{-4}\left(1-\frac {s^2}4\right)^2\left(1+\frac {s^2}4\right)d\omega_0dtds
\]
where $d\omega_0$ stands for the volume element for the round unit
sphere $S^2$
\begin{gather*}
 \text{vol}(\{ s > \epsilon\})  = 8\pi^2\int_\epsilon^2
s^{-4}\left(1-\frac {s^2}2 + \frac {s^4}{16}\right)\left(1+\frac {s^2}4\right)ds \\
\phantom{\text{vol}(\{ s > \epsilon\})}{}  = 8\pi^2\lambda \int_\epsilon^2 s^{-4}\left(1- \frac {s^2}4 - \frac
{s^4}{16} + \frac {s^6}{64}\right)ds \\
\phantom{\text{vol}(\{ s > \epsilon\})}{} = \frac {8\pi^2}3
\lambda\epsilon^{-3} - 2\pi^2 \lambda\epsilon^{-1} + 8\pi^2\lambda \left(
-\frac 1{3\times 8} + \frac 18 - \frac 18 + \frac 1{3\times 8}\right) +
O(\epsilon) \\
\phantom{\text{vol}(\{ s > \epsilon\})}{} = \frac {8\pi^2}3 \lambda \epsilon^{-3} - 2\pi^2
\lambda \epsilon^{-1} + O(\epsilon).
\end{gather*}
Thus
\[
V(S^1\times R^3, g_H^0) = 0.
\]

\noindent{\it c) AdS-Schwarzschild spaces}: \quad We
recall on $S^2\times R^2$
\[
g_{+1}^m = \left(1+r^2-\frac {2m}r\right)dt^2 + \frac {dr^2}{1+r^2-\frac {2m}r}
+ r^2(d\theta^2 + \sin^2\theta d\phi^2).
\]
First let us f\/ind the special def\/ining function, i.e.\ to have
\[
\frac 1{1+r^2 - \frac {2m}r}dr^2 = s^{-2}ds^2
\]
that is, if denote by $r = \rho/s$, where $\rho = \rho(s)$,
\[
\rho - s\rho' = \sqrt{\rho^2 + s^2 - 2ms^3/\rho},
\]
and $\rho(0) = 1$. One may solve it in power series
\[
\rho = 1 - \frac 14 s^2 + \frac m3 s^3 + \cdots.
\]
Then
\[
g^m_{+1} = s^{-2}\left(ds^2 + \left(\rho^2+s^2 - \frac {2ms^3}\rho\right)dt^2 +
\rho^2 (d\theta^2 + \sin^2\theta d\phi^2)\right).
\]
Note that $s\in [\epsilon, s_h]$ for $r\in [r_h, M_\epsilon]$,
\[
\log s_h = \log\epsilon + \int_{r_h}^{M_\epsilon} \frac 1
{\sqrt{1+r^2-\frac {2m}r}} dr < + \infty,
\]
and
\[
M_\epsilon  = \epsilon^{-1}\rho(\epsilon) = \epsilon^{-1}\left(1 - \frac
14 \epsilon^2 + \frac m3 s^3 + \cdots\right).
\]
Therefore
\begin{gather*}
 \text{vol}(\{ s > \epsilon \})  = \int_\epsilon^{s_h}
\int_{S^1(\lambda)} \int_{S^2} s^{-4} \sqrt{\rho^2 + s^2 - \frac
{2ms^2}\rho} \rho^2
dtd\sigma_0 ds \\
\phantom{\text{vol}(\{ s > \epsilon \})}{}  = 8\pi^2\lambda \int^{s_h}_\epsilon s^{-4} \sqrt{\rho^2 + s^2
-\frac {2ms^3}\rho} \rho^2 ds \\
\phantom{\text{vol}(\{ s > \epsilon \})}{} = 8\pi^2 \lambda \int_{r_h}^M s^{-1}\sqrt {1+r^2
-\frac {2m}r} r^2 \left(- \frac {ds}{dr}\right) dr \\
\phantom{\text{vol}(\{ s > \epsilon \})}{} = 8\pi^2\lambda\int_{r_h}^M r^2 dr = \frac {8\pi^2\lambda}3 (M^3-
r^3_h).
\end{gather*}
Thus the renormalized volume
\[
V(R^2 \times S^2, g^m_{+1}) = \frac {8\pi^2}3 \frac
{r^2_h(1-r^2_h)}{3r^2_h + 1},
\]
where $V(R^2 \times S^2, g^m_{+1}) < 0$ when $r_h > 1$; $V(R^2
\times S^2, g^m_{+1}) = 0$ only when $r_h = 1$ or $0$; and it
achieves its maximum value at $r_h = 1/\sqrt 3$
\[
V(R^2 \times S^2, g^m_{+1})_{\max} = \frac 19\cdot \frac {4\pi^2}3
\chi(R^2\times S^2).
\]

\noindent{\it d) AdS-Kerr spaces}: \quad We will omit
the calculation here. The renormalized volume
\[
V(\text{CP}^2\setminus \{p\}, g_\alpha)  = 4\pi^2 E_\alpha\left(-\frac 16
E_\alpha(\alpha^3 + 3\alpha^{-1}) + \frac 23 (\alpha +\alpha^{-1})
\right).
\]
Clearly, $V(\text{CP}^2\setminus \{p\}, g_\alpha)$ goes to zero when
$\alpha$ goes to $2$, and $V(\text{CP}^2\setminus \{p\}, g_\alpha)$
goes to $-\infty$ when $\alpha$ goes to $\infty$. One may f\/ind the
maximum value for the renormalized volume is achieved at $\alpha = 2
+ \sqrt 3$. Therefore
\[
V(\text{CP}^2\setminus \{p\}, g_\alpha)_{\max} = \frac {4\pi^2}3
\cdot \frac {2(4-\sqrt 3)}9 < \frac 12 \cdot \frac {4\pi^2}3
\chi(\text{CP}^2\setminus \{p\}).
\]

\subsection[Renormalized volume and Chern-Gauss-Bonnet formula]{Renormalized volume and Chern--Gauss--Bonnet formula}

We start with the Gauss--Bonnet formula on a surface $(M^2, g)$
\[
4\pi \chi(M^2) = \int_M K dv_g,
\]
where $K$ is the Gaussian curvature of $(M^2, g)$. The
transformation of the Gaussian curvature under a conformal change of
metrics $g_w = e^{2w}g$ is governed by the Laplacian as follows:
\[
- \Delta_g w + K[g] = K[e^{2w}g] e^{2w}.
\]
The Gauss--Bonnet formula for a compact surface with boundary $(M^2
g)$ is
\[
4\pi\chi(M) = \int_M K dV_g + 2\int_{\partial M} kd\sigma_g,
\]
where $k$ is the geodesic curvature for $\partial M$ in $(M, g)$.
The transformation of the geodesic curvature under a conformal
change of metric $g_w = e^{2w}g$ is
\[
-\partial_n w + k[g] = k[e^{2w}g] e^w,
\]
where $\partial_n$ is the inward normal derivative. Notice that
\[
 -\Delta [e^{2w} g]  = e^{-2w}(-\Delta [g]), \qquad
        - \partial_n [e^{2w}g]  = e^{-w} (-\partial_n[g]),
\]
for which we say they are conformally covariant. In four dimension
there is a rather complete analogue. We may write the
Chern--Gauss--Bonnet formula in the form
\[
8\pi^2 \chi(M^4) = \int_M (|W|^2 + Q)dV_g
\]
for closed 4-manifold and
\[
8\pi^2\chi(M^4) = \int_M(|W|^2 + Q)dV_g + 2\int_{\partial M} (L +
T)d\sigma_g,
\]
where $W$ is the Weyl curvature, $L$ is a point-wise conformal
invariant curvature of $\partial M$ in $(M, g)$.
\begin{gather*}
 Q  = \frac 16 (R^2 - 3|{\rm Ric}|^2 - \Delta R), \\
T  = - \frac 1{12} \partial n R + \frac 16 RH - R_{\alpha n \beta
n}L_{\alpha\beta} + \frac 19 H^3 - \frac 13 \text{Tr} L^3 - \frac 13
\tilde \Delta H,
\end{gather*}
$R$ is the scalar curvature, ${\rm Ric}$ is the Ricci curvature, $L$ is
the second fundamental form of $\partial M$ in $(M, g)$. We know the
transformation of $Q$ under a conformal change metric $g_w = e^{2w}
g$ is
\[
P_4[g] w + Q[g] = Q[e^{2w}g] e^{4w},
\]
where
\[
P_4 = (-\Delta)^2 + \delta\left\{\frac 23 Rg - 2\text{Ric} \right\}d
\]
is the so-called Paneitz operator, and the transformation of $T$ is
\[
P_3 [g] w + T[g] = T[e^{2w}]e^{3w},
\]
where
\begin{gather*}
P_3 = \frac 12 \partial_n \Delta_g - \tilde\Delta\partial_n + \frac
23 H \tilde\Delta + L_{\alpha\beta}\tilde\nabla_\alpha
\tilde\nabla_\beta + \frac 13 \tilde\nabla_\alpha H \cdot
\tilde\nabla_\alpha + \left(F-\frac 13 R\right)\partial_n.
\end{gather*}
We also have
\[ P_4[e^{2w} g]  = e^{-4w} P_4[g],\qquad
         P_3[e^{2w} g]  = e^{-3w}P_3[g].
\]

On the other hand, to calculate the renormalized volume in general,
for odd $n$, upon a choice of a special def\/ining function $x$, one
may solve
\[
-\Delta v = n \qquad\text{in $X^{n+1}$}
\]
for
\[
v = \log x + A + B x^n,
\]
$A$, $B$ are even in $x$, and $A|_{x=0}=0$. Let
\[
B_n [g, \hat g] = B|_{x=0}.
\]
Fef\/ferman and Graham observed

\begin{lemma}\label{lemma 2.5.1}
\[
V (X^{n+1}, g) = \int_M B_n [g, \hat g] dv[\hat g].
\]
\end{lemma}

We observe that the function $v$ in the above is also good in
conformal compacti\-f\/i\-cations. For example, given a conformally
compact Einstein 4-manifold $(X^4, g)$, let us consider the
compactif\/i\-cation $(X^4, e^{2v}g)$. Then
\[
Q_4 [e^{2v}g] = 0
\]
and its boundary is totally geodesic in $(X^4, e^{2v}g)$. Moreover
\[
T [e^{2v} g] = 3 B_3[g, \hat g].
\]
Therefore we obtain easily the following generalized
Chern--Gauss--Bonnet formula.

\begin{proposition}\label{proposition 2.5.2} Suppose that $(X^4, g)$ is a conformally
compact Einstein manifold. Then
\[
8\pi^2\chi(X^4) = \int_{X^4} (|W|^2 dv)[g]  +  6 V(X^4, g).
\]
\end{proposition}

\subsection{Topology of conformally compact Einstein
4-manifolds}

In the following let us summarize some of our works appeared in
\cite{CQY2}. From the generalized Chern--Gauss--Bonnet formula,
obviously
\[
V \leq \frac {4\pi^2}3 \chi(X)
\]
and the equality holds if and only if $(X^4, g)$ is hyperbolic.
Comparing with Chern--Gauss--Bonnet formula for a closed 4-manifold
\[
\frac 1{8\pi^2}\int_{M^4} (|W|^2 + \sigma_2)dv = \chi(M^4)
\]
one sees that the renormalized volume replaces the role of the
integral of  $\sigma_2$. In the following we will report some
results on the topology of a conformally compact Einstein 4-manifold
in terms of the size of the renormalized volume relative to the
Euler number, which is analogous to the  results of
Chang--Gursky--Yang \cite{CGY1, CGY2} on a closed 4-manifold with
positive scalar curvature and large integral of $\sigma_2$ relative
to the Euler number. The proofs mainly rely on the conformal
compactif\/ications discussed earlier, a simple doubling argument and
applications of the above mentioned results of Chang--Gursky--Yang
\cite{CGY1, CGY2}.

\begin{theorem}\label{theorem 2.6.1} Suppose $(X^4, g)$ is a conformally compact
Einstein $4$-manifold with its conformal infinity of positive Yamabe
constant and the renormalized volume $V$ is positive. Then $H^1(X,
R) = 0$.
\end{theorem}

\begin{theorem}\label{theorem 2.6.2} Suppose $(X^4, g)$ is a conformally compact
Einstein $4$-manifold with conformal infinity of positive Yamabe
constant. Then
\[
V > \frac 13  \frac {4\pi^2}3 \chi (X)
\]
implies that $H^2(X, R)$ vanishes.
\end{theorem}

A nice way to illustrate the above argument is the following. We may
consider the modif\/ied Yamabe constant
\[
Y^\lambda(M, [g]) = \inf_{g \in [g]} \frac {\int_M (R[g] + \lambda
|W^+|_g)dv_g}{ (\int_M dv_g)^{\frac {n-2}n}}.
\]
Then, one knows that $(M, [g])$ is of positive $Y^\lambda(M, [g])$
if and only if there is a metric $g \in [g]$ with $R + \lambda |W^+|
> 0$. As a consequence of the following Bochner formula
\[ \Delta \frac 12 |\omega|^2  = |\nabla \omega|^2 -
2W^+(\omega, \omega) + \frac 13 R|\omega|^2
  \geq |\nabla \omega|^2 + (R - 2\sqrt{6}|W^+|)|\omega|^2
\]
for any self-dual harmonic 2-form $\omega$, one easily sees that a
closed oriented 4-manifold with $Y^{-2\sqrt{6}} >0$ has its $b^+_2
=0$. We also observe

\begin{theorem}\label{theorem 2.6.3} Suppose $(X^4, g)$ is a conformally compact
Einstein $4$-manifold with its conformal infinity of positive Yamabe
constant and that
\[
V > \frac 12 \frac {4\pi^2}3 \chi(X).
\]
Then $X$ is diffeomorphic to $B^4$ and more interestingly $M$ is
diffeomorphic to $S^3$.
\end{theorem}

The detailed proofs of the above theorems are in our paper
\cite{CQY2}. One may recall
\begin{gather*} a)  \quad V(H^4, g_H) = \frac {4\pi^2}3,
\\
b) \quad V(S^1\times R^3, g_H) = 0,
\\
c) \quad V(S^2\times R^2, g_+^m) = \frac {8\pi^2}3 \frac
{r^2_h(1-r^2_h)}{3r^2_h + 1}  \leq \frac 19 \frac {4\pi^2}3\chi(S^2\times R^2), \\
d) \quad V(CP^2\setminus B, g_K) \leq \frac {4\pi^2}3
  \frac {2(4-\sqrt 3)}9   < \frac 13 \frac {4\pi^2}3\chi(CP^2\setminus B).
\end{gather*}
Theorem~\ref{theorem 2.6.2} is rather sharp, in cases (c) and (d)
the second homology is nontrivial while the renormalized volume is
very close to one-third of the maximum.

\pdfbookmark[1]{References}{ref}
\LastPageEnding


\begin{thebibliography}{99}

\footnotesize\itemsep=0pt

\bibitem{An} Anderson M., Orbifold compactness for spaces of Riemannian metrics
and applications,
{\it Math. Ann.} {\bf 331} (2005), 739--778,
\href{http://arxiv.org/abs/math.DG/0312111}{math.DG/0312111}.

\bibitem{AC} Anderson M., Cheeger J., Dif\/feomorphism f\/initeness for
manifolds with Ricci curvature and $L^{n/2}$-norm of curvature
bounded, {\it Geom. Funct. Anal.} {\bf 1} (1991),
231--252.

\bibitem{BC} Brezis  H.,  Coron J.M., Convergence of solutions of
H-systems or how to blow bubbles, {\it Arch. Ration. Mech. Anal.}
{\bf 89} (1985), 21--56.

\bibitem{BN} Bray  H., Neves A., Classif\/ication of prime
3-manifolds with Yamabe invariant larger than $RP^3$, {\it Ann. of
Math. (2)} {\bf 159} (2004), 407--424.

\bibitem{CQY1} Chang S.-Y.A., Qing J., Yang P., On a
conformal gap and f\/initeness theorem for a class of four-manifolds.
{\it Geom. Funct. Anal.} {\bf 17} (2007), 404--434, \href{http://arxiv.org/abs/math.DG/0508621}{math.DG/0508621}.

\bibitem{CQY2} Chang S.-Y.A., Qing  J., Yang P., On the
topology of conformally compact Einstein 4-manifolds, in
Noncompact Problems at the Intersection of Geometry, Analysis, and
Topology, {\it Contemp. Math.} {\bf 350} (2004),  49--61, \href{http://arxiv.org/abs/math.DG/0305085}{math.DG/0305085}.

\bibitem{CGY1} Chang S.-Y.A., Gursky M., Yang P.,  An
equation of Monge--Amp\'ere type in conformal geometry and 4-manifolds
of positive Ricci curvature, {\it Ann. of Math. (2)} {\bf 155} (2002),
709--787, \href{http://arxiv.org/abs/math.DG/0409583}{math.DG/0409583}.

\bibitem{CGY2} Chang S.-Y.A., Gursky M., Yang P., An apriori estimate for a fully nonlinear equation on 4-manifolds, {\it J.~D'Analyse Math.} {\bf 87} (2002), 151--186.


\bibitem{G} Graham C.R., Volume and area renormalizations for
conformally compact Einstein metrics, in The Proceedings of the 19th
Winter School ``Geometry and Physics'' (1999, Srn\`{i}), {\it Rend.
Circ. Mat. Palermo (2)}  {\bf 63} (2000), suppl., 31--42, \href{http://arxiv.org/abs/math.DG/9909042}{math.DG/9909042}.

\bibitem{GL} Graham  C.R., Lee J., Einstein metrics with
prescribed conformal inf\/inity on the ball, {\it Adv. Math.} {\bf 87}
(1991), 186--225.



\bibitem{Q} Qing J.,  On singularities of the heat f\/low for
harmonic maps from surfaces into spheres, {\it Comm. Anal. Geom.} {\bf 3} (1995),  297--315.

\bibitem{Q1} Qing J., On the rigidity for conformally compact
Einstein manifolds,  {\it Int. Math. Res. Not.} {\bf 2003} (2003), 1141--1153, \href{http://arxiv.org/abs/math.DG/0305084}{math.DG/0305084}.

\bibitem{St} Struwe M., Global compactness result for elliptic
boundary value problem involving limiting nonlinearities, {\it Math.
Z.} {\bf 187} (1984), 511--517.

\bibitem{TV-1} Tian  G., Viaclovsky J., Bach f\/lat
asymptotically ALE metrics, {\it Invent. Math.} {\bf 160} (2005),
357--415, \href{http://arxiv.org/abs/math.DG/0310302}{math.DG/0310302}.

\bibitem{TV-2} Tian  G., Viaclovsky J., Moduli space of
critical Riemannian metrics in dimension 4, {\it Adv. Math.} {\bf
196} (2005), 346--372, \href{http://arxiv.org/abs/math.DG/0312318}{math.DG/0312318}.


\end{thebibliography}
\end{document}